\documentclass[11pt]{article}
\usepackage{latexsym}

\usepackage[top=1in, bottom=1in, left=1in, right=1in]{geometry}

\usepackage{adjustbox}

\usepackage{amssymb}

\usepackage{graphicx} 
\usepackage{color} 
\usepackage{amsmath, amsthm, amssymb}
\usepackage{enumerate}
\usepackage{float}
\usepackage{url}
\usepackage{stackrel}
\usepackage{mathrsfs,dsfont}
\usepackage{caption}
\usepackage{subcaption}
\usepackage{wrapfig}
\usepackage{bbm}
\usepackage{bm}
\usepackage{epigraph}
\usepackage{physics}

\usepackage{mdwlist}

\usepackage{pifont}

\newtheorem{theorem}{Theorem}[section]
\newtheorem{corollary}[theorem]{Corollary}

\newtheorem{remark}[theorem]{Remark}

\newtheorem{examplex}[theorem]{Example}
\newenvironment{example}
  {\pushQED{\qed}\examplex}
  {\popQED\endexamplex}

\theoremstyle{definition}
\newtheorem{definition}[theorem]{Definition}
\theoremstyle{definition}

\newcommand{\been}{\begin{enumerate}}
\newcommand{\enen}{\end{enumerate}}
\newcommand{\beit}{\begin{itemize}}
\newcommand{\enit}{\end{itemize}}

\def\xrin{\xrightarrow{t \to \infty}}

\makeatletter
\newcommand{\xrightleftarrows}[2]{%
  \mathrel{\mathop{%
    \vcenter{\offinterlineskip\m@th
      \ialign{\hfil##\hfil\cr
        \hphantom{$\scriptstyle\mspace{8mu}{#1}\mspace{8mu}$}\cr
        \rightarrowfill\cr
        \vrule height0pt width 2em\cr
        \leftarrowfill\cr
        \hphantom{$\scriptstyle\mspace{8mu}{#2}\mspace{8mu}$}\cr
        \noalign{\kern-0.3ex}
      }%
    }%
  }\limits^{#1}_{#2}}%
}
\makeatother

\def\DD{\mathcal{D}}

\def\spn{{\rm span}}

\def\ve{\varepsilon}

\def\ds{\displaystyle}

\newcommand{\R}{\mathbb{R}}

\usepackage{tikz}    
\usetikzlibrary{shapes,automata,positioning,arrows,fit}
\usepackage{mathtools}  
\usetikzlibrary{snakes}

\usepackage{caption}
\usepackage{tikz-cd}
\usetikzlibrary{babel}

\numberwithin{equation}{section}

\usepackage{tikz-3dplot}

\usepackage{relsize}

 \tikzset{every node/.style={auto}}
 \tikzset{every state/.style={rectangle, minimum size=0pt, draw=none, font=\normalsize}}
  \tikzset{bend angle=7}
  
  \usepackage{array}

\usepackage{authblk}

\begin{document}

\title{Power-engine-load form for \\ dynamic absolute concentration robustness}

\author{Badal Joshi \footnote{Department of Mathematics, California State University San Marcos (bjoshi$@$csusm.edu)} ~ and Gheorghe Craciun \footnote{Departments of Mathematics and Biomolecular Chemistry, University of Wisconsin-Madison (craciun$@$wisc.edu)}}

\date{}

\maketitle

\begin{abstract}
In a reaction network, the concentration of a species with the property of dynamic absolute concentration robustness (dynamic ACR) converges to the same value independent of the overall initial values. 
This property endows a biochemical network with output robustness and therefore is essential for its functioning in a highly variable environment. 
It is important to identify structure of the dynamical system as well as constraints required for dynamic ACR. 
We propose a power-engine-load form of dynamic ACR and obtain results regarding convergence to the ACR value based on this form. 
~

\noindent {\bf AMS subject codes:} 34C20, 37N25, 37N35, 92C42

\end{abstract}

\section{Introduction}

{\em Dynamic} absolute concentration robustness (dynamic ACR), introduced in \cite{joshi2022foundations}, is a property of dynamical systems wherein one variable converges to a unique value independent of the initial values. 
This variable is the {\em dynamic ACR variable} and the unique value is its {\em ACR value}. 
Dynamic ACR is significant for applications to biochemistry. 
Biochemical systems need to perform robustly in a wide variety of conditions. 
Dynamic ACR provides a mechanism for such robustness. 
In signal response circuits, an essential design feature is that the response depends on the signal but not on the internal state of the signaling circuit. 
Different initial states (which model the internal state of the signaling circuit when the signal arrives) generically result in different final states, especially when some conservation laws hold. 
However, if one of the variables is invariant across all final states, such as is the case in dynamic ACR, then this variable can be taken to be the invariant signal response. 
The dependence of the signal strength is encoded in the reaction rate constants, and therefore signal response will depend on the signal strength via the rate constants. 
However, the signal response will not depend on the initial values, i.e. the internal state of the signaling circuit.

Formally, dynamic ACR is the property that all trajectories (with some minor restrictions on allowed initial values) converge to the hyperplane $\{x_i = a_i^*\}$. 
It is desirable to identify sufficient conditions, either on the underlying reaction network or the dynamical system, that generate dynamic ACR. 
Necessary and sufficient conditions for dynamic ACR in complex balanced systems have been obtained in  \cite{joshi2022foundations}. 
Moreover, a classification of small networks with dynamic ACR appears in \cite{joshi2023motifs}. 

Shinar and Feinberg \cite{shinar2010structural} defined {\em static} absolute concentration robustness (static ACR) as the property that all positive steady states are located in the hyperplane $\{x_i = a_i^*\}$. 
Furthermore, they gave sufficient network conditions for the reaction network to have static ACR in some species. 
Additionally, they gave many biochemically realistic networks and showed that these networks satisfy their network conditions and therefore have static ACR. 
Many of the biochemical networks appearing in their paper may have the property of dynamic ACR as well, although this remains an open question. 
It is difficult to characterize dynamic ACR because it requires understanding the limiting behavior for arbitrary initial values. 
We propose a {\em power-engine-load form} (see \eqref{eq:normalformACR}) for the differential equation satisfied by the supposed dynamic ACR variable, which may help establish dynamic ACR in some situations as we show in this paper. 
In Theorem \ref{thm:nonzeroload}, we establish sufficient conditions for convergence to the ACR value, based on properties of the terms appearing in the power-engine-load form equation. 

Dynamic ACR has been experimentally observed in bacterial two-component signaling systems such as the EnvZ-OmpR system and the IDHKP-IDH system \cite{russo1993essential,hsing1998mutations,batchelor2003robustness,shinar2007input,alon2019introduction,hart2013utility}. 
Absolute concentration robustness (both static and dynamic) is related to the concept called {\em robust perfect adaptation} \cite{kim2020absolutely,xiao2018robust,cappelletti2020hidden} studied from a control theory perspective. 
For a biochemical perspective on perfect adaptation, see \cite{khammash2021perfect, sontag2003adaptation}. 
Structural requirements for robust perfect adaptation in biomolecular networks are studied in \cite{gupta2022universal}. 
While our goal in this paper is not to explore the connections between the two fields of dynamical systems and control theory, we mention the above papers here as a help to any readers interested in understanding the connection further. 

This article is organized as follows.  
Section~\ref{sec:reactionnetworks}  gives the background on reaction networks, mass action kinetics, and dynamic absolute concentration robustness (ACR). 
Section~\ref{powerengineloadform} introduces a power-engine-load form for dynamic ACR, examples, and convergence results. The main result of this paper is Theorem \ref{thm:nonzeroload}.
Section~\ref{sec:applications} studies the dynamics of several reaction networks where we apply power-engine-load form and illustrate the use of Theorem  \ref{thm:nonzeroload}.
In many cases, the reaction network is coupled with inflows. While each network requires special analysis, we ultimately apply Theorem \ref{thm:nonzeroload} to show dynamic ACR for the network under study.

\section{Reaction networks} \label{sec:reactionnetworks}

An example of a {\em reaction} is $A + B \to 2C$, which is a schematic representation of the process where a molecule of species $A$ and a molecule of species $B$ react with one another and result in two molecules of a species $C$. The abstract linear combination of species $A+B$ that appears on the left of the reaction arrow is called the {\em reactant complex} while $2C$ is called the {\em product complex} of the reaction $A+B \to 2C$. 
We assume that for any reaction the product complex is different from the reactant complex. 
A {\em reaction network} is a nonempty, finite set of species and a nonempty, finite set of reactions such that every species appears in at least one complex.  
We will use {\em mass action kinetics} for the rate of reactions. 

\subsection{Mass action systems}

In mass action kinetics, each reaction occurs at a rate proportional to the product of the concentrations of species appearing in the reactant complex. 
We conventionally use lower case letters $a, b, c$ to denote the species concentrations of the corresponding species $A, B, C$, respectively. 
Under mass action kinetics, the reaction $A+B \to 2C$ occurs at rate $k ab$ where $k$ is the reaction rate constant, conventionally placed near the reaction arrow, as follows: $A+B \xrightarrow{k} 2C$. 
Consider the reaction network $\{ A+ B \xrightarrow{k_1} 2C, ~2C \xrightarrow{k_2} 2A \}$. 
Application of mass action kinetics results in the following dynamical system, called {\em mass action system}, that describes the time-evolution of the species concentrations. 
\begin{align*}
\dot a &= -k_1 ab + 2 k_2 c^2 \\
\dot b &= -k_1 ab  \\
\dot c &= 2k_1 ab - 2 k_2 c^2 
\end{align*}
For further details on reaction networks and mass action systems, see for instance \cite{joshi2015survey,yu2018mathematical}. We use standard notation and terminology of dynamical systems, such as steady states, stability, basin of attraction etc., see for instance \cite{david2018ordinary, guckenheimer2013nonlinear}.

\subsection{Dynamic absolute concentration robustness} 

Shinar and Feinberg defined (static) absolute concentration robustness (static ACR) as the property that all positive steady states are contained in a fixed hyperplane parallel to a coordinate hyperplane \cite{shinar2010structural}.  
This means that for every positive steady state, one of the coordinates is invariant. 
Several network/structural conditions for static ACR have been identified \cite{shinar2010structural, meshkat2022absolute, joshi2022foundations, joshi2023motifs}. 
Static ACR is designed to model the property of a robust signal response despite variability in the internal state of the signaling circuit. 
However, the property accurately models output robustness only if there is convergence to a positive steady state for every initial value, and indeed, such convergence is not required in the definition of static ACR. 

We defined dynamic absolute concentration robustness (dynamic ACR) as the property that all initial values lead to convergence to a fixed hyperplane parallel to a coordinate hyperplane \cite{joshi2022foundations} -- a requirement somewhat more general than converging to a steady state in the hyperplane. We believe that this property better captures the idea of output robustness, when compared to static ACR, since it models long-term dynamics.  
Network conditions for dynamic ACR in small reaction networks were found in \cite{joshi2023motifs}. 
The condition found there is geometric in nature -- for dynamic ACR in networks with two reactions and two species, the reactant polytope (line segment joining the two reactant complexes in their geometric embedding) must be parallel to a coordinate axis. 
We seek an extension of this condition in arbitrary networks. If the differential equation for the candidate ACR species has a certain special form (power-engine-load form), then it may have dynamic ACR provided that the power and load satisfy some additional conditions.

Dynamic ACR was defined for autonomous systems in \cite{joshi2022foundations}.  
We generalize to include non-autonomous systems.
Throughout the paper, we assume that $\DD$ is a dynamical system defined by $\dot x = F(x,t)$ with $x \in \R^n_{\ge 0}$ for which $\R^n_{\ge 0}$ is forward invariant. 
\begin{definition} \label{def:compatible}
The {\em kinetic subspace} of $\DD$ is defined to be the linear span of the image of $F$, denoted by $\spn(\Im(F))$. 
Two points $x,y \in \R^n_{\ge 0}$ are {\em compatible} if $y - x \in \spn(\Im(F))$. 
The sets $S, S' \subseteq \R^n_{\ge 0}$ are {\em compatible} if there are $x \in S$ and $x' \in S'$ such that $x$ and $x'$ are compatible.
A {\em compatibility class $S$} is a nonempty subset of $\R^n_{\ge 0}$ such that $x,y \in S$ if and only if $y - x \in \spn(\Im(F))$.
\end{definition}
\begin{definition} \label{def:dyn_acr}
$\DD$ is a {\em dynamic ACR system} if there is an $i \in \{1,\ldots, n\}$ with $F_i \not \equiv 0$ and a positive $a_i^* \in \R_{> 0}$ such that for any $(t_0,x(t_0)) \in \R \times \R^n_{> 0}$ with $x(t_0)$ is compatible with $\{x \in \R^n_{>0} ~|~ x_i = a_i^*\}$, a unique solution to $\dot x = F(x,t)$ exists up to some maximal $T_0(t_0,x(t_0)) \in (t_0, \infty]$, and $x_i(t) \xrightarrow{t \to T_0} a_i^*$. Any such $x_i$ and $a_i^*$  is a {\em dynamic ACR variable} and its {\em dynamic ACR value}, respectively.  
\end{definition}
If the dynamical system $\dot x = F(x,t)$ does not have the possibility of a finite-time blow-up, then $T_0(t_0,x(t_0)) = \infty$ for any $(t_0,x(t_0)) \in \R \times \R^n_{>0}$. In this paper, we focus on dynamic ACR in systems where a unique solution is assumed to exist for all positive time.

\section{Power-engine-load form of dynamic ACR} \label{powerengineloadform}

In \cite{joshi2023motifs}, we discuss the network conditions for static and dynamic ACR in reaction networks with two species and two reactions. One of the conditions is that the reactant polytope (i.e. in the Euclidean embedding of the reaction network, the line segment joining the two reactant complexes) is parallel to one of the coordinate axes. A natural generalization of this geometric condition is the one given below, called  {\em power-engine-load form for dynamic ACR}. 

Suppose that $x \in \R^n$ and for some $1 \le i \le n$, $x_i$ satisfies
\begin{align}\label{eq:normalformACR}
\dv{x_i}{t} ~=~ {\color{olive} \underbrace{f(x(t),t)}_{\mbox{power}}}~\cdot~ {\color{magenta} \underbrace{(x_i^*-x_i(t))}_{\mbox{engine}}} ~+~ {\color{cyan} \underbrace{g(x(t),t)}_{\mbox{load}}}
\end{align}
then under some reasonable conditions on ``power'' $f(x(t),t)$ and ``load'' $g(x(t),t)$, the variable $x_i$ will have dynamic ACR with the value $x_i^*$. 
An example of sufficient conditions is $g \equiv 0$, $f > 0$, and $\int_0^\infty f(x(t),t)~dt = \infty$. We state the precise result and prove other convergence results in Section \ref{sec:convergence}. 

\subsection{Examples of power-engine-load form in mass action systems}
\begin{example}
The simplest example of a mass action system with power-engine-load form is the reaction network 
$
0 \stackrel[k']{k}{\rightleftarrows} X
$
and its associated ODE 
\[
\dv{x}{t} = {\color{olive} k'} {\color{magenta} \left(k/k' - x \right)}. 
\]
It is straightforward to show that $x$ has dynamic ACR with value $k/k'$. 
\end{example}
\begin{example} \label{ex:archetypal}
Biologically interesting cases require a reactant complex with more than one species (so that a {\em `reaction'} can occur) and some positive mass conservation law involving all species.
The following well-known network is the simplest that satisfies these requirements of species interaction and mass conservation:
\begin{align} \label{eq:canonicalACRnet}
A+B \xrightarrow{k_1} 2B, \quad 
B \xrightarrow{k_2} A. 
\end{align}
The positive mass conservation law $a+b = const$ is apparent from the mass action ODE system:
\begin{align} \label{eq:canonicalACR_odes}
\dot a = {\color{olive} k_1b} {\color{magenta} \left(k_2/k_1 - a \right)}, \quad \dot b = - k_1b\left(k_2/k_1 - a \right). 
\end{align}
The ODE for $a$ has power-engine-load form and it is known (see for instance \cite{joshi2022foundations}) that $a(t) \xrin k_2/k_1$ for any initial value $(a(0), b(0)) = (a_0,b_0)$ that satisfies $a_0 + b_0 \ge k_2/k_1$.  
\end{example}
\begin{example}
In bacterial two-component signaling systems, the circuit mechanism for robust signal transduction from the cell environment to its interior involves a bifunctional component, a mechanism that is found in thousands of biological systems \cite{alon2019introduction}. 
One such system is the E. coli IDHKP-IDH glyoxylate bypass regulation system whose core ACR module is 
\begin{align} \label{eq:idhkpidhnet}
X+E & \stackrel[k_2]{k_1}{\rightleftarrows} C_1 \xrightarrow{k_3} Y+E \nonumber\\ 
Y+C_1 & \stackrel[k_5]{k_4}{\rightleftarrows} C_2 \xrightarrow{k_6} X+C_1. 
\end{align}
It is known that \eqref{eq:idhkpidhnet} has static ACR in $Y$ \cite{shinar2010structural}. 
The mass action ODE equation for the concentration of $Y$ can be written in power-engine-load form, using the fact that the static ACR value of $y$ is $k_3/k_4(1+k_5/k_6)$, as follows. 
\begin{align*}
\dot y &= k_3c_1-k_4c_1y+k_5c_2 \\
&= k_3c_1-k_4c_1y+k_5\frac{k_3}{k_6}c_1 + \frac{k_5}{k_6}\left (k_6 c_2 - k_3c_1 \right) \\
&= {\color{olive} k_4 c_1} {\color{magenta} \left( \frac{k_3}{k_4}\left( 1+\frac{k_5}{k_6}\right) - y\right)} + {\color{cyan} \frac{k_5}{k_6}\left (k_6 c_2 - k_3c_1 \right)}. 
\end{align*}
The load is not identically zero in this case. 
We will prove in future work that $y$ has dynamic ACR with ACR value $k^* =  \dfrac{k_3}{k_4}\left( 1+\dfrac{k_5}{k_6}\right)$.
As shown in Theorem \ref{thm:nonzeroload}, a sufficient condition for dynamic ACR is  that for any positive initial value $\left (k_6 c_2(t) - k_3c_1(t) \right) \xrin 0$  and that $\int_0^\infty c_1(t)~dt = \infty$. 
Checking that the two conditions are satisfied requires some additional ideas which will be developed in future work.
\end{example}

We now discuss conditions on ``power'' $f(x(t),t)$ and ``load'' $g(x(t),t)$ that either ensure or prevent dynamic ACR. 

\subsection{Convergence results for power-engine-load form dynamic ACR} \label{sec:convergence}

\begin{theorem}[Zero load] \label{thm:zeroload}
Consider the dynamical system $\dot x = F(x,t)$ with continuously differentiable
 $F: \R^n_{\ge 0} \times \R_{\ge 0} \to \R^n$ and for which 
 $\R^n_{\ge 0}$ is forward invariant.
Suppose there is an $x_i^* \in \R_{> 0}$ such that $F_i(x,t)|_{x \in \R^n_{>0}} = 0$ if and only if $x_i = x_i^*$.  
The following hold for every solution $(x(t))_{t \ge 0}$ of $\dot x = F(x,t)$.
\been
\item $x_i^* - x_i(t)$ has the same sign as $x_i^* - x_i(0)$ for all $t \ge 0$.  
\item Suppose $x_i(0) \ne x_i^*$. Then $\ds \frac{F_i(x(t),t)}{x_i^* - x_i(t)}$ has the same sign as $\ds \frac{F_i(x(0),0)}{x_i^* - x_i(0)}$ for all $t \ge 0$.  
\item Suppose $x_i(0) \ne x_i^*$. 
\been
\item If $\ds \frac{F_i(x(0),0)}{x_i^* - x_i(0)} > 0$ then $\abs{x_i(t) - x_i^*}$ is strictly decreasing on $[0,\infty)$ and 
\[
\abs{x_i(\infty) - x_i^*} = \abs{x_i(0) - x_i^*} \exp\left(-\int_0^\infty \frac{F_i(x(s),s)}{x_i^*-x_i(s)} ds\right).
\] 
\item If $\ds \frac{F_i(x(0),0)}{x_i^* - x_i(0)} < 0$ then $\abs{x_i(t) - x_i^*}$ is strictly increasing on $[0,\infty)$. 
\enen
\item Suppose $x_i(0) \ne x_i^*$. Then $x_i(t) \xrin x_i^*$ if and only if 
\[
\int_0^\infty \frac{F_i(x(t),t)}{x_i^* - x_i(t)}  dt = \infty.
\]
\enen
\end{theorem}
\begin{proof}
If $x_i(0) = x_i^*$ then $\dot x_i|_{t =0} = F_i(x(0),0) = 0$ and so $x_i(t) = x_i^*$ for all $t \ge 0$.   
If there is a $t>0$ such that $x_i(t) - x_i^*$ and $x_i(0) - x_i^*$ have different signs then by continuity of $x_i(t)$, there must be a $t' \in (0,t)$ such that $x_i(t') = x_i^*$. But then $x_i(t) = x_i^*$ for all $t \in \R$, which is a contradiction. 

If $x_i(0) \ne x_i^*$ then by the previous part $\ds \frac{F_i(x(t),t)}{x_i^* - x_i(t)}$ is defined for all time $t \ge 0$. Since $F_i(x,t) \ne 0$ for $x_i \ne x_i^*$, the second result follows.  
Since $x_i(t) \ne x_i^*$ for all $t \ge 0$, we can divide by $x_i^* - x_i$ and integrate to get  
\begin{align}
\frac{dx_i}{x_i^*-x_i} = \frac{F_i(x(t),t)}{x_i^*-x_i(t)}dt 
\implies \frac{x_i^*-x_i(t)}{x_i^*-x_i(0)} = \exp\left(-\int_0^t \frac{F_i(x(s),s)}{x_i^*-x_i(s)} ds\right). \label{eq:intontheright}
\end{align}
By the previous result, the integrand has the same sign as $\ds \frac{F_i(x(0),0)}{x_i^* - x_i(0)}$ for all positive time, and so result 3 follows. Finally, result 4 follows from taking $\lim_{t \to \infty}$ on both sides. 
\end{proof}

\begin{example}
Consider the mass action system \eqref{eq:canonicalACR_odes}, $\dot a = k_1b \left(k - a \right), ~ \dot b = - k_1b\left(k - a \right)$. The variable $a$ has power-engine-load form  with zero load and power $k_1b > 0$ for every $(a,b) \in \R^2_{> 0}$.  
It is clear that $a(0) = k$ if and only if $\dot a|_{t\ge0} = 0$. So assume that $a(0) \ne k$. 

The mass action system with initial value $(a(0), b(0)) = (a_0, b_0) \in \R^2_{\ge 0}$ can be solved explicitly after using the conservation relation $a(t) + b(t) = a_0 + b_0$, thus reducing to the one-dimensional system $\dot b = -k_1b (k + b - a_0 - b_0)$. We get
\begin{align*}
b(t) = \begin{cases} 
\dfrac{a_0 + b_0 -k}{1 + \left(\frac{a_0 - k}{b_0} \right) e^{-(a_0 + b_0 -k)t}} &\mbox{ if } a_0 + b_0 \ne k, \\ 
\dfrac{k-a_0}{1+ (k-a_0) t} &\mbox{ if } a_0 + b_0 = k.
\end{cases}
\end{align*}
and $a(t) = a_0 + b_0 - b(t)$.  

For $(a_0,b_0) \in \R^2_{> 0}$, it is well-known that $a(t) \xrin k$ if and only if $a_0 + b_0 \ge k$. 
We now argue that $\int_0^\infty b(t)~dt = \infty$ if and only if $a_0 + b_0 \ge k$ and $a_0 \ne k$. 
Indeed, for any $a_0 < k$, 
\[
\int_0^\infty \frac{k-a_0}{1+ (k-a_0) t}~dt = \infty. 
\]
Moreover
\begin{align} \label{eq:ohgeognogn}
\lim_{t \to \infty} \frac{a_0 + b_0 -k}{1 + \left(\frac{a_0 - k}{b_0} \right) e^{-(a_0 + b_0 -k)t}} = 
\begin{cases}
a_0 + b_0 -k &\mbox{ if } a_0 + b_0 > k, \\
0 &\mbox{ if } a_0 + b_0 < k,
\end{cases}
\end{align}
implies that
\[
\int_0^\infty \frac{a_0 + b_0 -k}{1 + \left(\frac{a_0 - k}{b_0} \right) e^{-(a_0 + b_0 -k)t}}~dt  
\begin{cases}
= \infty &\mbox{ if } a_0 + b_0 > k, \\
< \infty &\mbox{ if } a_0 + b_0 < k,
\end{cases}
\]
because in the first case ($a_0 + b_0 > k$) the integrand does not go to zero, the integral is clearly divergent and in the latter case ($a_0 + b_0 < k$) the integrand is $\approx \exp{(a_0+b_0-k)t}$, so the integral is convergent. 

This shows that $a \xrin k$ if and only if either $a(0) = k$ or $\int_0^\infty k_1 b(t)~dt = \infty$. 
\end{example}
\begin{corollary} \label{cor:zeroload}
Suppose the hypotheses of Theorem \ref{thm:zeroload} hold and $x_i(0) \ne x_i^*$. 
If 
\begin{align} \label{eq:limitcondition}
\liminf_{t \to \infty} \frac{t F_i(x(t),t)}{x_i^*-x_i(t)} \in (0,\infty], 
\end{align}
then $x_i(t) \xrin x_i^*$. 
\end{corollary}
\begin{proof}
By \eqref{eq:limitcondition}, there is a $\delta > 0$ and a $t_0 \ge 0$ such that $\dfrac{tF_i(x(t),t)}{x_i^*-x_i(t)} > \delta$ for all $t \ge t_0$.  Then for any $T \ge t_0$, 
\begin{align*}
\int_{0}^\infty \frac{F_i(x(t),t)}{x_i^* - x_i}  dt 
\ge \int_{t_0}^T \frac{F_i(x(t),t)}{x_i^* - x_i}  dt  
= \int_{t_0}^T \frac{1}{t}~\frac{t F_i(x(t),t)}{x_i^* - x_i}  dt  
> \delta  \int_{t_0}^T \frac{1}{t} dt \xrightarrow{T \to \infty} \infty, 
\end{align*}
and so by Theorem \ref{thm:zeroload}, $x_i \to x_i^*$. 
\end{proof}
\begin{example}
The condition \eqref{eq:limitcondition} is sufficient but not necessary for $x_i \to x_i^*$. Consider the mass action system of the following reaction network:
\begin{align}
2X_1 \xrightarrow{k_1} X_1, &\quad 2X_2 \xrightarrow{k_2} X_2,  \nonumber\\
X_1 + X_2 \xrightarrow{k_2} X_2,  &\quad Y + X_1 \stackrel[k_4]{k_5}{\rightleftarrows} X_1, 
\end{align}
where the rate constants are the same for the reactions $2X_2 \to X_2$ and $X_1 + X_2 \to X_2$ ($X_2$ degrades both $X_1$ and $X_2$ at the same rate). The mass action ODEs are
\begin{align}
\dot x_2 &= -k_2x_2^2, \nonumber\\
\dot x_1 &= -k_1x_1^2 - k_2x_1x_2, \nonumber\\
\dot y &= k_5x_1\left(\frac{k_4}{k_5} - y\right). 
\end{align}
It is simple to check that given an arbitrary positive initial value $(y(0), x_1(0), x_2(0)) = (b_0,b_1,b_2) \in \R^3_{> 0}$, the unique solution satisfies 
\begin{align*}
x_2(t) &= \frac{b_2}{1+k_2b_2t}, \\
x_1(t) &= \frac{b_1b_2k_2}{(1+b_2k_2t)(b_2k_2+b_1k_1\log(1+b_2k_2t)} \sim \frac{1}{k_1t\log(1+b_2k_2t)} 
\end{align*}
and so $\int_{t_0}^\infty x_1(t) dt= \infty$ for any $t_0 > 0$ which implies that $y \to k_4/k_5$ by Theorem \ref{thm:zeroload}. However, 
$\liminf_{t \to \infty}  (t x_1(t)) = 0$, so Corollary \ref{cor:zeroload} does not apply. 
\end{example}


When the load $g$ is nonzero, the power $f$ must overpower the load $g$ for convergence of $x_i$ to $x_i^*$.

\begin{theorem} \label{thm:nonzeroload}
Consider the dynamical system $\dot x = F(x,t)$ with 
 $F: \R^n_{\ge 0} \times \R_{\ge 0} \to \R^n$ and for which 
 $\R^n_{\ge 0}$ is forward invariant.
Suppose that for some $i \in \{1,\ldots, n\}$ we have
\[
F_i(x,t) = f(x,t) \cdot \left( x_i^* - x_i \right) + g(x,t),
\] 
with $f(x,t) > 0$ in $\R^n_{> 0} \times \R_{> 0}$, and 
$g \not \equiv 0$. 
Let $(x(t))_{t \ge 0}$ be a solution of $\dot x = F(x,t)$ such that 
 $\ds \int_0^\infty f(x(t),t) dt = \infty$, and such that the limit 
$\ds \alpha \coloneqq  \lim_{t \to \infty} \frac{g(x(t),t)}{f(x(t),t)}$ exists, with $\alpha > - x_i^*$. 
Then we have $x_i(t) \xrin x_i^* + \alpha$.
\end{theorem}
\begin{proof}
Note that the equation  $\ds \dv{x_i}{t}  = F_i(x,t)$ can be rewritten as
\begin{align} \label{eq:givensystemodes_ver_2}
\dv{x_i}{t} &= f(x,t) \left(x_i^* + \frac{g(x,t)}{f(x,t)} - x_i\right).
\end{align}
Then, if we denote $\widetilde x_i^* = x_i^* + \alpha$ and $\widetilde g = g/f - \alpha$, we can reduce our problem to showing that if $x_i(t)$ satisfies the equation 
\begin{align} \label{eq:givensystemodes_ver_3}
\dv{x_i}{t} &= f(x,t) \left(\ \widetilde x_i^*  - x_i + \widetilde g(x,t) \ \right),
\end{align}
and $\ds  \lim_{t \to \infty} \widetilde g(x(t),t) = 0$, then $x_i(t) \xrin \widetilde x_i^*$.

\bigskip

For any fixed $\varepsilon \in (0, \widetilde x_i^*)$ we will now show that there exists some $T_0>0$ such that $x_i(t) \in (\widetilde x_i^* - \varepsilon, \widetilde x_i^* + \varepsilon)$ for all $t > T_0$.
For this, let us first choose some $T_1 > 0$ such that  $\widetilde g(x(t),t)  \in ( - \frac{\varepsilon}{2},  \frac{\varepsilon}{2})$ for all $t > T_1$. Note that this implies that the interval $(\widetilde x_i^* - \varepsilon,  \widetilde x_i^* + \varepsilon)$ is an invariant set of   \eqref{eq:givensystemodes_ver_2} for $t>T_1$. Therefore, if $x_i(T_1) \in (\widetilde x_i^* - \varepsilon,  \widetilde x_i^* + \varepsilon)$, we can choose $T_0 = T_1$. Assume now that $x_i(T_1) \notin (\widetilde x_i^* - \varepsilon,  \widetilde x_i^* + \varepsilon)$, and for example $x_i(T_1) \ge \widetilde x_i^* + \varepsilon$ (the case where $x_i(T_1) \le \widetilde x_i^* - \varepsilon$ is analogous).

Let us assume  that the inequality $x_i(t) \ge \widetilde x_i^* + \varepsilon$ holds for all $t > T_1$; we will show that this leads to a contradiction. Indeed, for any such $t$ we have 
$$
\widetilde x_i^*  - x_i + \widetilde g(x,t) < -\frac{\varepsilon}{2},
$$
which implies that 
\begin{align} 
\dv{x_i}{t} &\le -\frac{\varepsilon}{2} f(x,t),
\end{align}
for all $t > T_1$; but this, together with the hypothesis that $\ds \int_0^\infty f(x(t),t) dt = \infty$, would imply that $\ds \lim_{t \to \infty}x_i(t) = -\infty$, which contradicts our assumption that for all $t > T_1$ we have $x_i(t) \ge \widetilde x_i^* + \varepsilon$. 

Therefore we obtain the desired conclusion that $x_i(t)$ will enter the interval $(\widetilde x_i^* - \varepsilon, \widetilde x_i^* + \varepsilon)$ in some finite time $T_0$, and will remain inside it for all $t > T_0$. 
\end{proof}

\section{Applications of power-engine-load form} \label{sec:applications}

We use the power-engine-load form and Theorem \ref{thm:nonzeroload} to prove dynamic ACR in certain reaction networks taken with inflows. 
The underlying reaction networks have dynamic ACR when there are no inflows or outflows, i.e. when the system is closed. 
We show that even when certain inflows are included, the property of dynamic ACR persists. Moreover, the ACR value for the open system (with inflows) remains the same in many instances as the ACR value for the closed/isolated system. 
In Section \ref{subsec:outflowstoo}, we also consider an application where both inflows and outflows are included, and again given the right conditions the dynamic ACR property persists. 
In Section \ref{subsec:bifun}, we consider a simple enzyme catalysis network with a bifunctional enzyme. 
We show once again that dynamic ACR persists under many different possible inflow conditions. 
Each individual network requires special analysis in combination with application of Theorem \ref{thm:nonzeroload}. 

The reaction networks that we consider here have a conserved quantity when the system is closed. 
When inflows (but not outflows) are included, the total concentration grows in an unbounded manner. 
This behavior can be realized in a chemostat (for some finite time) but there are also biologically realistic situations where this occurs. For instance, specific ion channels on the surface of a cell might be open and selectively permeable which results in influx of channel-specific ions from outside the cell volume. 
Our analysis shows that {\em there are invariant quantities even while influxes are ongoing}, that is even when the overall system is not near any long-term equilibrium. 
One variable, the ACR variable, can consistently converge to a robust value even when the overall trajectory does not converge. 
Robust convergence of one variable when the overall system is far from equilibrium can even occur when inflows and outflows are both present,  as shown in Section \ref{subsec:outflowstoo}. 
In all these cases, Theorem \ref{thm:nonzeroload} is used to show convergence.

\subsection{Motifs of static and dynamic ACR} \label{subsec:motifs}

\begin{figure}[h!] 
\centering
\begin{tikzpicture}[scale=1]

\draw [->, line width=1.25] ({-6+0.5},{0}) -- ({-6-0.5},{0-1});
\draw [->, line width=1.25] ({-6-0.5},{0}) -- ({-6+0.5},{0+1});
\draw [-, line width=1, dashed] ({-6-0.5},{0}) -- ({-6+0.5},{0});

\draw [->, line width=1.25] ({-3+0.5},{0}) -- ({-3-0.5},{0});
\draw [->, line width=1.25] ({-3-0.5},{0}) -- ({-3+0.5},{0});
\draw [-, line width=1, dashed] ({-3-0.5},{0}) -- ({-3+0.5},{0});

\draw [->, line width=1.25] ({0+0.5},{0}) -- ({0-0.5},{0+1});
\draw [->, line width=1.25] ({0-0.5},{0}) -- ({0+0.5},{0-1});
\draw [-, line width=1, dashed] ({0-0.5},{0}) -- ({0+0.5},{0});

\end{tikzpicture}
\caption{Motifs with two reactions and two species that have both static and dynamic ACR in one species $A$.}
\label{fig:all3motifs}
\end{figure}
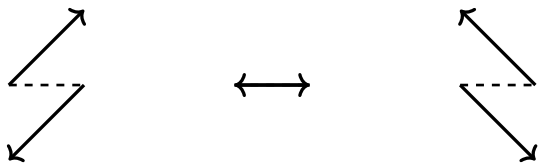

In \cite{joshi2023motifs}, we identified the network motifs with two reactions and two species which have both static ACR and dynamic ACR when treated as isolated or closed systems. 
Each motif represents an entire infinite set of reaction networks. 
We consider three such motifs shown in Figure \ref{fig:all3motifs} and we select one network for each motif.  It is shown in \cite{joshi2023motifs} that these three are the only motifs (up to interchange of species labels) which have two reactions, two species, and have both properties of static ACR and dynamic ACR. 

In each network motif, the arrows represent reactions while the (dashed) line segment joining the two arrow tails is the so-called reactant polytope which plays an important role in the classification of ACR systems \cite{joshi2023motifs} -- the reactant polytope is parallel to a coordinate axis in all ACR systems with two reactions and two species. 

Each static and dynamic ACR motif (with two reactions and two species) has the following properties: the reactant polytope is parallel to a coordinate axis, the two reaction arrows are mutually parallel and they point towards each other. The three motifs are then characterized by the slope of the reaction vectors, which can be positive, zero, or negative. 
We will refer to the motifs accordingly as ``positive/zero/negative slope motif''. 

For each motif, we choose a particular network and show that when considered as an open system by allowing inflows, the ACR property is preserved in many cases. 
One representative of negative slope motif, see \eqref{eq:canonicalACRnet}, is often studied as a simple example of a conservative system with static ACR. 
We show that not only does it have dynamic ACR, it has many other striking robustness properties.

\subsection{Positive slope motif}
Our first example is a network based on the motif ~
\begin{tikzpicture}[scale=0.5]
\draw [->, line width=1.25] ({-6+0.5},{0}) -- ({-6-0.5},{0-1});
\draw [->, line width=1.25] ({-6-0.5},{0}) -- ({-6+0.5},{0+1});
\draw [-, line width=1, dashed] ({-6-0.5},{0}) -- ({-6+0.5},{0});
\end{tikzpicture}.
~
For this motif, we choose a particular reaction network whose true/non-flow reactions are $\{A+B \to 0, B \to A+2B\}$, and we include the inflows $\{0 \to A, 0 \to B\}$ with arbitrary time-dependent rates. 
So the reaction network is 
\begin{align*}
A+B \xrightarrow{k_1} 0, \quad &
B \xrightarrow{k_2} A+2B, \\
0    \xrightarrow{g_a(t)} A,   \quad &
0    \xrightarrow{g_{b}(t)} B. 
\end{align*}
whose mass action system is 
\begin{equation} \label{eq:motif1}
\begin{aligned} 
\dot a &=  k_1b(k^*-a) + g_a(t),\\
\dot b &= k_1b(k^*-a) + g_b(t), 
\end{aligned}
\end{equation}
where $k^* = k_2/k_1$.

\begin{theorem}
Consider the mass action system \eqref{eq:motif1} with time-dependent inflow rates $g_a: \R_{\ge 0} \to \R_{\ge 0}$ and $g_b: \R_{\ge 0} \to \R_{\ge 0}$, such that $g_a(t)$ is bounded and 
\begin{equation} \label{eq:intcond1}
\int_{0}^{\infty} \left(g_b(t) - g_a(t)\right) dt  = \infty,
\end{equation}
then $a(t) \xrin k^*$. 
\end{theorem}
\begin{proof}
Since $\dot b - \dot a = g_b - g_a$, by \eqref{eq:intcond1}, we have that $b(t) - a(t) \xrin \infty$, which implies that $b(t) \xrin \infty$. Then the conditions of Theorem \ref{thm:nonzeroload} hold for $a$ and from 
\[
\lim_{t \to \infty} \frac{g_a(t)}{b(t)} = 0, 
\]
we conclude that $a(t) \xrin k^*$. 
\end{proof}

\begin{corollary}
Consider the mass action system \eqref{eq:motif1} with constant inflows $g_b > g_a$. Then $a(t) \xrin k^*$. 
\end{corollary}

\subsection{Negative slope motif}

A representative of the motif ~
\begin{tikzpicture}[scale=0.5]
\draw [->, line width=1.25] ({0+0.5},{0}) -- ({0-0.5},{0+1});
\draw [->, line width=1.25] ({0-0.5},{0}) -- ({0+0.5},{0-1});
\draw [-, line width=1, dashed] ({0-0.5},{0}) -- ({0+0.5},{0});
\end{tikzpicture}
~
was used to illustrate static ACR by Shinar and Feinberg \cite{shinar2010structural}. 
We show here that this network has even stronger robustness properties than previously known. 
In particular, even in face of inflows which send the total concentration to infinity, the species with ACR still converges to a finite value, and moreover, the finite value is the same as its ACR value when there are no inflows. 

Consider the following open reaction network: 
\begin{equation}
\begin{aligned} 
A+B \xrightarrow{k_1} 2B, \quad &
B \xrightarrow{k_2} A, \\
0    \xrightarrow{g_a(t)} A,   \quad &
0    \xrightarrow{g_{b}(t)} B, 
\end{aligned}
\end{equation}
whose mass action system is 
\begin{equation} \label{eq:motif2}
\begin{aligned} 
\dot a &=  k_1b(k^*-a) + g_a(t),\\
\dot b &= -k_1b(k^*-a) + g_b(t), 
\end{aligned}
\end{equation}
where $k^* = k_2/k_1$.

If $g_a$ and $g_b$ are constant, then the ACR value survives. We give a more general result below, allowing $g_a$ and $g_b$ to be arbitrary functions of time. 
One implication of the result is that the ACR value survives for arbitrary functions, provided that $g_a(t)$ does not increase too fast. 

\begin{theorem} \label{thm:flowmotif2exp}
Consider the mass action system \eqref{eq:motif2} with time-dependent inflow rates $g_a: \R_{\ge 0} \to \R_{\ge 0}$ and $g_b: \R_{\ge 0} \to \R_{\ge 0}$, such that 
\begin{equation} \label{eq:intcond3}
G(t) \coloneqq \int_{0}^{t} \left(g_a(s) + g_b(s)\right) ds  \xrin \infty,
\end{equation}
and 
\begin{equation} \label{eq:flowcond3}
g_a(t)/G(t) \xrin \alpha.  
\end{equation}
then $a(t) \xrin k^* + \alpha/k_1$ for any $(a(0),b(0)) \in \R^2_{> 0}$. 
\end{theorem}
\begin{proof}
We will first show that $b/G \to 1$. Note that $\dot a + \dot b = g_a + g_b$ implies that 
$
a(t) + b(t) = a(0) + b(0) + G(t), 
$
and so $a+b \xrin \infty$. 

Define the following invertible change of coordinates: 
\begin{equation} \label{eq:changeofcoods}
\begin{aligned}
\R_{\ge 0}^2 \setminus \{(0,0)\} &\to \R_{>0} \times [0,1] \\
(a,b) &\mapsto \left(x ,\beta \right) = \left(a+b, b/(a+b)\right)
\end{aligned}
\end{equation}
In $(x,\beta)$ coordinates, the dynamical system \eqref{eq:motif2} is: 
\begin{equation}
\begin{aligned} 
\dot x &= g_a(t) + g_b(t), \\
\dot \beta &= k_1\beta \left(x (1-\beta) - k^* \right) + \frac{(1-\beta)g_b - \beta g_a}{x}. 
\end{aligned}
\end{equation}
%
%
%

The first equation has the solution $x(t) = x(0) + G(t)$, which implies that $x(t)$ grows monotonically to infinity.
It then follows that for any $\delta \in (0,\frac{1}{2})$ there is a time $T_\delta>0$ such that for all $t>T_\delta$ we have $\dot \beta(t) > 1$ whenever $\beta(t) \in (\delta, 1-\delta)$, because under such assumptions  $x(t)$ is large enough to imply that the positive term 
\[
k_1\beta \left(x (1-\beta)\right)
\]
is much larger than all the negative terms on the right-hand side of the equation for $\dot \beta(t)$. 
Moreover, we obtain that if $\beta(0) \in  (\delta, 1-\delta)$ then the solution  $\beta(t)$ becomes larger than $1 - \delta$ for all $t > T_\delta+1$, which implies that $\beta(t) \to 1$, and therefore $b/G \to 1$. 

Finally, we check the hypotheses of Theorem \ref{thm:nonzeroload}. 
We may rewrite the equation for $\dot a$ in \eqref{eq:motif2} as 
\begin{equation} 
\begin{aligned} 
\dot a &=  k_1b\left(\frac{\alpha}{k_1} + k^*-a \right) + G(t) \left(\frac{g_a(t)}{G(t)} - \alpha \frac{b}{G(t)} \right). 
\end{aligned}
\end{equation}
Clearly $b \to \infty$ and so $\int b \to \infty$. Moreover, we have that
\[
\lim_{t \to \infty} \frac{G(t)}{b(t)} \left(\frac{g_a(t)}{G(t)} - \alpha \frac{b}{G(t)} \right) = \lim_{t \to \infty} \left(\frac{g_a(t)}{G(t)} - \alpha  \right) = 0. 
\]
Therefore we may conclude from Theorem \ref{thm:nonzeroload} that 
$a(t) \xrin k^* + \alpha/k_1$ for any $(a(0),b(0)) \in \R^2_{> 0}$. 
\end{proof}

\begin{remark}
Theorem \ref{thm:flowmotif2exp} is formally about a two-dimensional non-autonomous system. 
We show in examples below that the theorem can be applied to study the dynamics of higher dimensional systems. 
See Examples \ref{ex:polygrowth}, \ref{ex:expgrowth} and \ref{ex:tetrationgrowth}. 
\end{remark}

\begin{corollary} \label{cor:polygrowth}
Consider the mass action system \eqref{eq:motif2} with time-dependent inflow rates $g_a: \R_{\ge 0} \to \R_{\ge 0}$ and $g_b: \R_{\ge 0} \to \R_{\ge 0}$, such that $g_a(t)$ is a polynomial of $t$ (that takes only positive values). 
This includes the case of constant inflows. 
Then $a(t) \xrin k^*$ for any $(a(0),b(0)) \in \R^2_{> 0}$. 
\end{corollary}
\begin{proof}
We only need to check the hypotheses of Theorem \ref{thm:flowmotif2exp}. Clearly \eqref{eq:intcond3}  holds. 
Moreover for \eqref{eq:flowcond3}, we have 
\begin{align*}
0 \le \frac{g_a(t)}{G(t)} \le \frac{g_a(t)}{\int_0^t g_a(s) ds} \xrin 0, 
\end{align*}
and so $\alpha = 0$. It follows that $a \xrin k^*$. 
\end{proof}

\begin{example} \label{ex:polygrowth}
Consider the following reaction network for some positive integer $d$ 
\begin{equation} \label{net:polygrowth}
\begin{aligned} 
0 &\to X_1,\\
X_j &\to X_j + X_{j+1}, \quad (1 \le j \le d-1)\\
X_d + Y &\xrightarrow{k_1} 2Y, \\
Y &\xrightarrow{k_2} X_d. 
\end{aligned}
\end{equation}
Note that $d=1$ is allowed and corresponds to constant inflow of $X_d = X_1$. 
For $j \in \{1, \ldots, d-1\}$, $X_j \sim t^j$ and so Corollary \ref{cor:polygrowth} implies that $X_d \to k_2/k_1$. 
\end{example}

Now we consider the case where the inflows are exponentially growing in time. 
\begin{example} \label{ex:expgrowth}
Consider the following reaction network
\begin{equation} \label{net:expgrowth}
\begin{aligned} 
Z \xrightarrow{\alpha} 2Z, &\quad Z \to Z + A\\
A + B &\xrightarrow{k_1} 2B, \\
B &\xrightarrow{k_2} A. 
\end{aligned}
\end{equation}
The $Z$ could represent an exponentially growing virus which produces a toxin $A$. The differential equation for $Z$ is $\dot z = \alpha z$ which has the solution $z(t) =  e^{\alpha t}$ (assuming that $z(0) = 1$) and so the concentrations of $A$ and $B$ evolve according to:
\begin{equation}
\begin{aligned} 
\dot a &= k_1 b (k^* - a) + e^{\alpha t}, \\
\dot b &= -k_1 b (k^* - a). 
\end{aligned}
\end{equation}
Now 
\[
G(t) = \int_0^t  e^{\alpha s} ds =   \frac{e^{\alpha t} - 1}{\alpha}, 
\]
which clearly goes to infinity and moreover $g_a(t)/G(t) \to \alpha$. 
By Theorem \ref{thm:flowmotif2exp}, it then follows that 
\[
a(t) \xrin k^* + \alpha/k_1.
\]
A similar result holds even when both $A$ and $B$ have exponentially growing inflows. For instance, for the system 
\begin{equation}
\begin{aligned} 
\dot a &= k_1 b (k^* - a) + e^{\alpha t}, \\
\dot b &= -k_1 b (k^* - a) + e^{\beta t}, 
\end{aligned}
\end{equation}
we have that 
\begin{align*}
\lim_{t \to \infty} a(t) = \begin{cases}
k^* + \alpha/k_1 &\mbox{ if } \alpha > \beta \\
k^* + \alpha/(2k_1) &\mbox{ if } \alpha = \beta \\
k^*  &\mbox{ if } \alpha < \beta. 
\end{cases}
\end{align*}
\end{example}

\begin{remark}
It is possible to find an explicit solution for the case of arbitrary time-dependent inflow in $A$, $g_a(t)$ and no inflow in $B$. We give such a solution here. 
One can see from the form of the solution, even though it is explicit, that the limiting value of $a(t)$ given in Theorem \ref{thm:flowmotif2exp} is not obtained easily. 
Moreover, we emphasize the additional shortcoming that the explicit solution is for the restricted situation of no inflow in $B$. 

Consider the following system: 
\begin{align}
\dot a = -k_1 ab + k_2 b + g_a(t),  \quad 
\dot b = k_1 ab - k_2 b. 
\end{align}
Let $G_a(t) = \int_0^t g_a(s) ds$, so $\dot a + \dot b = g(t)$ has the solution $a(t) + b(t) = a(0) + b(0)  + G_a(t)$. We can use this to write an ODE in $a$ only: 
\begin{align}
\dot a = -k_1(a - k^*)(a(0) + b(0) + G_a(t) - a ) + g_a(t), 
\end{align}
where $k^* = k_2/k_1$.

\begin{theorem}
Let $g_a(t)$ be an arbitrary function of time $t$ and $G_a(t) = \int_0^t g_a(s)~ds$. 
For any initial value $a(0) = a_0 \ge 0$ and $b(0) = b_0 \ge 0$, 
\begin{align} \label{eq:odeequation}
\dot a = -k_1(a - k)(a_0 + b_0 + G_a(t) - a ) + g_a(t)
\end{align}
has the solution
\begin{align} \label{eq:odesolution}
a(t) = a_0 + b_0 + G_a(t) - \frac{b_0 q(t)}{1+k_1 b_0 Q(t)}
\end{align}
for all time $t \ge 0$ where 
\begin{align}
q(t) = \exp\left[ k_1(a_0+b_0-k)t + k_1\int_0^t G_a(s)~ds \right], 
\end{align}
and $Q(t) = \int_0^t q(s)~ds$. 
\end{theorem}
\begin{proof}
The proof is a simple verification. From their definitions, $q(0) = 1, G_a(0) = 0, Q(0) = 0$ and so 
\[
a_0 + b_0 + G_a(0) - \frac{b_0 q(0)}{1+k_1 b_0 Q(0)} = a_0 + b_0 - b_0 = a_0, 
\]
so the initial condition is satisfied. 
We note that $\dot q(t) = k_1 q(t) \left( a_0 + b_0 - k + G_a(t) \right)$ and so 
\begin{align*}
\dot a(t) 
&= g_a(t) - b_0 \frac{k_1 q(t) \left( a_0 + b_0 - k + G_a(t) \right)}{1+k_1 b_0 Q(t)} + k_1 \left(\frac{b_0 q(t)}{1+k_1b_0 Q(t)} \right)^2, 
\end{align*}
and therefore
\begin{align*}
&\dot a + k_1(a - k)(a_0 + b_0 + G_a(t) - a ) - g_a(t) \\
&= \dot a + k_1\left( a_0 + b_0 + G_a(t) - \frac{b_0 q(t)}{1+k_1 b_0 Q(t)} - k\right) \left( \frac{b_0 q(t)}{1+k_1 b_0 Q(t)} \right) - g_a(t) = 0,
\end{align*}
which completes the verification. 
\end{proof}

\end{remark}

\subsection{Negative slope motif with inflows and outflows} \label{subsec:outflowstoo}

We can extend the results in the previous section even to the situation where one or both species are in outflow. 
We consider the case where both $A$ and $B$ are in outflow and at equal rates. 
Outflows with equal rates for all species is a reasonable (and standard) assumption for chemostats. 
We will allow the inflows to be arbitrary functions of time. 
Suppose we have the following reaction network: 
\begin{align*}
0    \xrightarrow{g_a(t)} A,   \quad &
0    \xrightarrow{g_{b}(t)} B,  \\
A    \xrightarrow{\ell} 0, \quad & 
B    \xrightarrow{\ell} 0,  \\
A+B \xrightarrow{k_1} 2B, \quad &
B \xrightarrow{k_2} A, 
\end{align*}
whose mass action system is: 
\begin{equation} \label{eq:equaloutflows}
\begin{aligned}
\dot a &= k_1b(k^* - a) + g_a(t) - \ell a, \\ 
\dot b &= -k_1b(k^* - a) + g_b(t) - \ell b. 
\end{aligned}
\end{equation}

\begin{figure}[h!]
\centering
\includegraphics[scale=0.41]{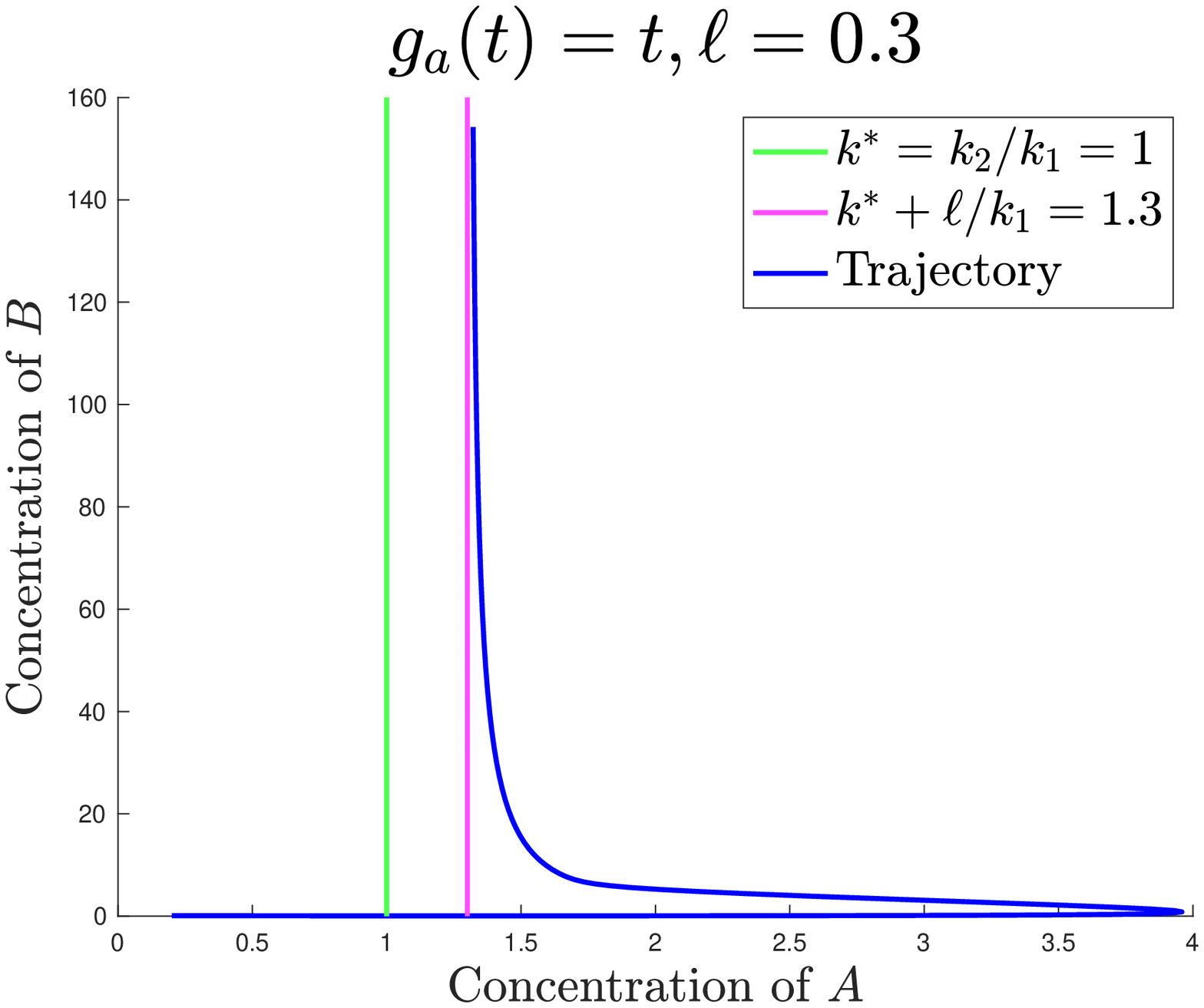}
\includegraphics[scale=0.41]{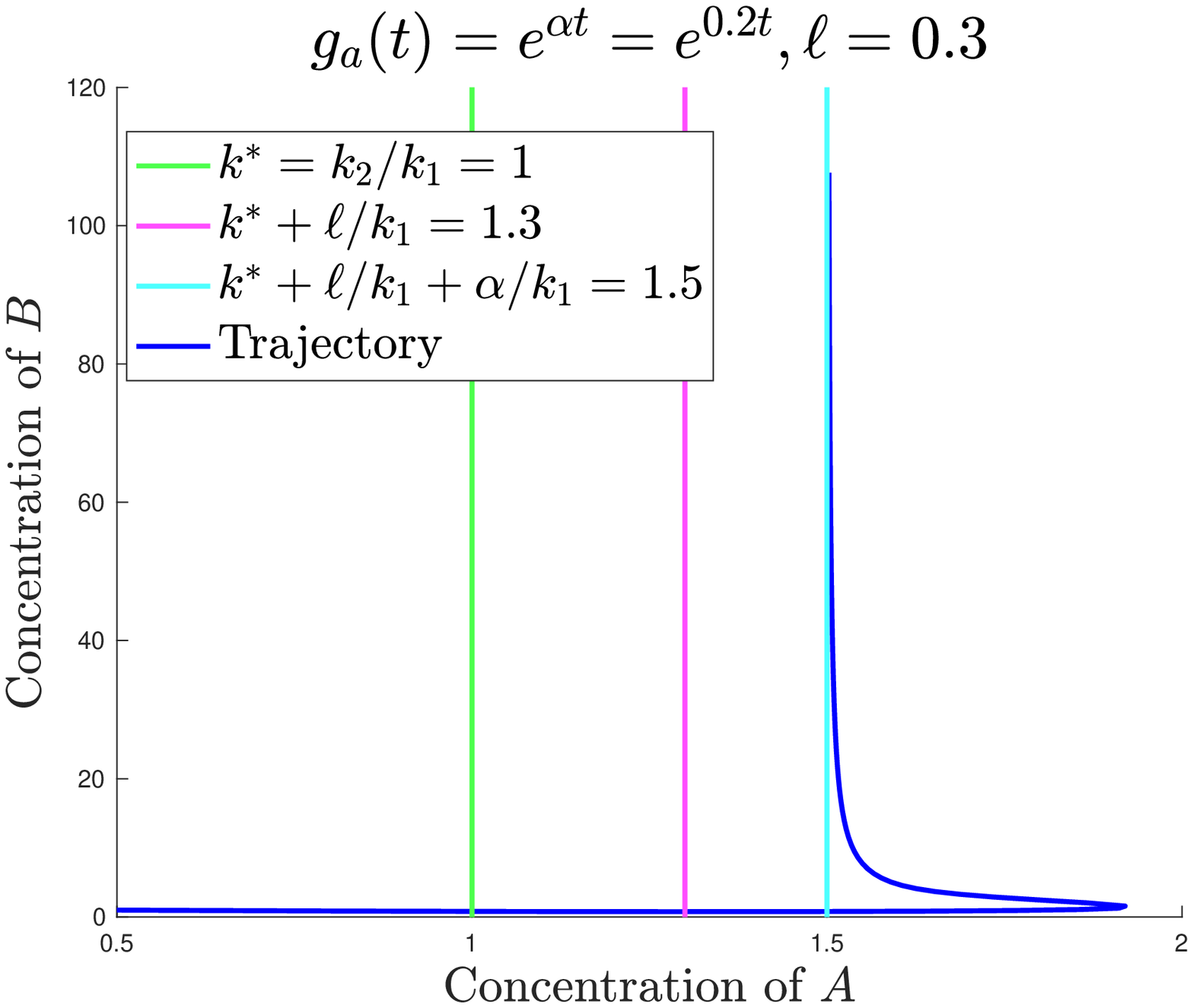}
\caption{The reaction network $\{A+B \xrightarrow{k_1} 2B, ~ B \xrightarrow{k_2} A\}$ coupled with an inflow $0 \xrightarrow{g_a(t)} A$ (with $g_a(t) \xrin \infty$) and outflows $A \xrightarrow{\ell} 0, ~ B \xrightarrow{\ell} 0$.
When the inflow $g_a(t)$ is sub-exponential, $a$ converges to a value that is equivalent to changing the rate constant of the reaction $B \to A$ from $k_2$ to $k_2 + \ell$. When the inflow $g_a(t) = \exp(\alpha t)$ is exponential, $a$ converges to a value that is equivalent to changing the rate constant of the reaction $B \to A$ from $k_2$ to $k_2 + \alpha + \ell$.}
\label{fig:bothoutflows}
\end{figure}

\begin{theorem} \label{thm:equaloutflows}
Consider \eqref{eq:equaloutflows} with $g_a, g_b:\R_{\ge 0} \to \R_{\ge 0}$ such that $g(t) \coloneqq g_a(t)+g_b(t) \xrin \infty$, and $\ell > 0$. 
Assume that if $g_b \equiv 0$ then $b(0) > 0$. 
Let 
\[
H(t)  \coloneqq  e^{-\ell t} \int_0^t g(s) e^{\ell s} ~ds.
\] 
If $\alpha  \coloneqq  \lim_{t \to \infty}  g_a(t)/H(t)$ exists and is finite, then 
\[
a(t) \xrin k^* + \frac{\alpha}{k_1}. 
\]
\end{theorem}
\begin{proof}
We will first show that $b/H \to 1$. 
Define the following invertible change of coordinates: 
\begin{align*}
\R_{\ge 0}^2 \setminus \{(0,0)\} &\to \R_{>0} \times [0,1] \\
(a,b) &\mapsto \left(x ,\beta \right) = \left(a+b, b/(a+b)\right). 
\end{align*}
In $(x,\beta)$ coordinates, the dynamical system \eqref{eq:equaloutflows} is
\begin{equation} \label{eq:equaloutflowsbeta}
\begin{aligned} 
\dot x &= g(t) - \ell x \nonumber\\
\dot \beta &= k_1 \beta \left(x (1-\beta) - k^* \right) + \frac{(1-\beta)g_b - \beta g_a}{x}
\end{aligned}
\end{equation}
The equation for $\dot x$ can be solved explicitly using an integrating factor: 
\[
x(t) = e^{-\ell t} \left( \int_0^t g(s) e^{\ell s} ~ds + x(0) \right) = x(0) e^{-\ell t} + H(t). 
\]
We now show that $x \to \infty$. By hypothesis, we have that $g \xrin \infty$. So for any $N>0$, there is a $T_N> (\ln 2)/ \ell$ such that $g(t) > N$ for $t > T_N$. So, 
\begin{align*}
x(2T_N) &\ge H(2T_N) \ge e^{-2 \ell T_N} \int_{T_N}^{2T_N} g(s) e^{\ell s}~ds \\ 
&\ge e^{-2 \ell T_N} N \int_{T_N}^{2T_N} e^{\ell s}~ds = N \frac{1- \exp(-\ell T_N)}{\ell} \ge \frac{N}{2 \ell}. 
\end{align*}
Since $N$ is arbitrary, $x \xrin \infty$. 

Now note that the equation for $\dot \beta$ is the {\em same} as the corresponding equation in Theorem \ref{thm:flowmotif2exp} (in particular, the outflow parameter $\ell$ does not appear in $\dot \beta$). 
Therefore, the same reasoning as in Theorem \ref{thm:flowmotif2exp} gives us that $b/x \xrin 1$ and $b/H \xrin 1$. 
Moreover, $a/x \to 0$ and 
\[
\frac{a}{b} = \frac{a/x}{b/x} \to 0. 
\]

Finally, we check the hypotheses of Theorem \ref{thm:nonzeroload}. 
Clearly $b \to \infty$ and so $\int b \to \infty$. Moreover, we have that
\[
\lim_{t \to \infty} \frac{- \ell a +  g_a(t)}{k_1b(t)} = \lim_{t \to \infty}\frac{1}{k_1} \left(- \ell \frac{a}{b} +  \frac{g_a(t)/H(t)}{b/H(t)} \right) = \frac{\alpha}{k_1}. 
\]
Therefore, we conclude from Theorem \ref{thm:nonzeroload} that 
$a(t) \xrin k^* + \alpha/k_1$ for any $b(0) > 0$ or $g_b \not \equiv 0$. 
\end{proof}
\begin{remark}
In Theorem \ref{thm:equaloutflows}, we assumed equal outflows of species $A$ and $B$, $\ell_a = \ell_b = \ell$, a standard assumption for a chemostat. Simulations show that even when this assumption is relaxed and the outflow rates are allowed to be different, the property of dynamic ACR in $A$ is preserved. 
We conjecture that the ACR value is obtained by replacing $\ell$ with $\ell_b$, the outflow of species $B$, in the expression for $H(t)$ in Theorem \ref{thm:equaloutflows}. 
\end{remark}
\begin{corollary} \label{cor:equaloutflows}
Assume the hypotheses and notation of Theorem \ref{thm:equaloutflows}. Suppose that $g_a$ is differentiable. Then 
\[
\alpha = \lim_{t \to \infty} \left[\ell \frac{g_a(t)}{g(t)} + \frac{g_a'(t)}{g(t)} \right], 
\]
and so 
\[
a(t) \xrin k^* + \frac{1}{k_1}\lim_{t \to \infty} \left[\ell \frac{g_a(t)}{g(t)} +  \frac{g_a'(t)}{g(t)}\right]. 
\]
\end{corollary}
\begin{proof}
Note that 
\begin{align*}
\alpha = \lim_{t \to \infty} \frac{g_a(t)}{H(t)} = \lim_{t \to \infty}  \frac{g_a(t) e^{\ell t}}{\int_0^t g(s) e^{\ell s}ds}.
\end{align*}
For any $g_a(t) \not \equiv 0$, both the numerator and the denominator go to infinity, so the result then follows from applying L'Hospital's rule. 
\end{proof}

\begin{example}
Consider the reaction network \eqref{net:polygrowth} in Example \ref{ex:polygrowth} with $d \ge 2$ (so that $g_a(t) \to \infty$), and additionally include outflows 
$X_d \xrightarrow{\ell} 0, ~ Y \xrightarrow{\ell} 0,$
so the complete reaction network is 
\begin{equation} \label{net:polygrowthoutflow}
\begin{aligned} 
0 &\to X_1,\\
X_j &\to X_j + X_{j+1}, \quad (1 \le j \le d-1)\\
X_d + Y &\xrightarrow{k_1} 2Y, \\
Y &\xrightarrow{k_2} X_d, \\
X_d& \xrightarrow{\ell} 0, \quad Y \xrightarrow{\ell} 0. 
\end{aligned}
\end{equation}
As earlier, for $j \in \{1, \ldots, d-1\}$, $X_j \sim t^j$. Choose rate constants and initial values so that $X_j = t^j$ for $j \in \{1, \ldots, d-1\}$. 
So the concentrations for $X_d$ and $Y$ follow: 
\begin{align*}
\dot x_d &= k_1y(k^* - x_d) + t^{d-1} - \ell x_d, \\
\dot y &= -k_1y(k^* - x_d)  - \ell y. 
\end{align*}
In the notation of Theorem \ref{thm:equaloutflows}, $g_{x_d}(t) = g(t) = t^{d-1}$ and by Corollary \ref{cor:equaloutflows}
\begin{align*}
\alpha = \ell + \lim_{t \to \infty} \frac{g_a'(t)}{g(t)} = \ell, 
\end{align*}
and so 
\[
x_d \xrin k^* + \frac{\ell}{k_1}. 
\]
See left panel of Figure \ref{fig:bothoutflows} for a representative trajectory. 
\end{example}

\begin{example}
Consider the reaction network \eqref{net:expgrowth} in Example \ref{ex:expgrowth} with outflows included, so the complete reaction network is 
\begin{equation}
\begin{aligned} 
Z \xrightarrow{\alpha} 2Z, &\quad Z \to Z + A\\
A + B &\xrightarrow{k_1} 2B, \\
B &\xrightarrow{k_2} A, \\
B \xrightarrow{\ell} 0, &\quad A \xrightarrow{\ell} 0. 
\end{aligned}
\end{equation}
The concentrations of $A$ and $B$ evolve according to:
\begin{equation}
\begin{aligned} 
\dot a &= k_1 b (k^* - a) + e^{\alpha t} - \ell a, \\
\dot b &= -k_1 b (k^* - a) - \ell b. 
\end{aligned}
\end{equation}
Noting that $g_a(t) = g(t) = e^{\alpha t}$, by Corollary \ref{cor:equaloutflows}, 
\[
a \xrin k^* + \frac{\ell}{k_1} + \frac{\alpha}{k_1}.  
\]
See right panel of Figure \ref{fig:bothoutflows} for a representative trajectory. 
\end{example}

\subsection{A representative of the negative slope motif with stronger robustness properties}

Another representative of the negative slope motif has even stronger robustness properties. 
Consider the following open reaction network: 
\begin{equation}
\begin{aligned} 
A+2B \xrightarrow{k_1} 3B, \quad &
2B \xrightarrow{k_2} A + B, \\
0    \xrightarrow{g_a(t)} A,   \quad &
0    \xrightarrow{g_{b}(t)} B, 
\end{aligned}
\end{equation}
whose mass action system is 
\begin{equation} \label{eq:motif5}
\begin{aligned} 
\dot a &=  k_1b^2(k^*-a) + g_a(t),\\
\dot b &= -k_1b^2(k^*-a) + g_b(t), 
\end{aligned}
\end{equation}
where $k^* = k_2/k_1$.

\begin{theorem} \label{thm:flowmotif5exp}
Consider the mass action system \eqref{eq:motif5} with time-dependent inflow rates $g_a: \R_{\ge 0} \to \R_{\ge 0}$ and $g_b: \R_{\ge 0} \to \R_{\ge 0}$, such that 
\begin{equation} \label{eq:intcond5}
G(t) \coloneqq \int_{0}^{t} \left(g_a(s) + g_b(s)\right) ds  \xrin \infty,
\end{equation}
and 
\begin{equation} \label{eq:flowcond5}
g_a(t)/(G(t))^2 \xrin \alpha.  
\end{equation}
Then $a(t) \xrin k^* + \alpha/k_1$ for any $(a(0),b(0)) \in \R^2_{> 0}$. 
\end{theorem}
\begin{proof}
The proof is analogous to that of Theorem \ref{thm:flowmotif2exp}. 
We will first show that $b/G \to 1$. Note that $\dot a + \dot b = g_a + g_b$ implies that 
$
a(t) + b(t) = a(0) + b(0) + G(t), 
$
and so $a+b \xrin \infty$. 
After changing to $(x,\beta)$ coordinates defined in \eqref{eq:changeofcoods}, the dynamical system \eqref{eq:motif2} is: 
\begin{equation}
\begin{aligned} 
\dot x &= g_a(t) + g_b(t), \\
\dot \beta &= k_1\beta^2 x \left(x (1-\beta) - k^* \right) + \frac{(1-\beta)g_b - \beta g_a}{x}. 
\end{aligned}
\end{equation}

The first equation has the solution $x(t) = x(0) + G(t)$, which implies that $x(t)$ grows monotonically to infinity.
It then follows that for any $\delta \in (0,\frac{1}{2})$ there is a time $T_\delta>0$ such that for all $t>T_\delta$ we have $\dot \beta(t) > 1$ whenever $\beta(t) \in (\delta, 1-\delta)$, because under such assumptions  $x(t)$ is large enough to imply that the positive term 
\[
k_1\beta^2 x \left(x (1-\beta)\right)
\]
is much larger than all the negative terms on the right-hand side of the equation for $\dot \beta(t)$. 
Moreover, we obtain that if $\beta(0) \in  (\delta, 1-\delta)$ then the solution  $\beta(t)$ becomes larger than $1 - \delta$ for all $t > T_\delta+1$, which implies that $\beta(t) \to 1$, and therefore $b/G \to 1$. 

Finally, we check the hypotheses of Theorem \ref{thm:nonzeroload}. 
We may rewrite the equation for $\dot a$ in \eqref{eq:motif2} as 
\begin{equation} 
\begin{aligned} 
\dot a &=  k_1b^2 \left(\frac{\alpha}{k_1} + k^*-a \right) + G^2(t) \left(\frac{g_a(t)}{G^2(t)} - \alpha \frac{b^2}{G^2(t)} \right). 
\end{aligned}
\end{equation}
Clearly $b \to \infty$ and so $\int b \to \infty$. Moreover, we have that
\[
\lim_{t \to \infty} \frac{G^2(t)}{b^2(t)} \left(\frac{g_a(t)}{G^2(t)} - \alpha \frac{b^2}{G^2(t)} \right) = \lim_{t \to \infty} \left(\frac{g_a(t)}{G^2(t)} - \alpha  \right) = 0. 
\]
Therefore we may conclude from Theorem \ref{thm:nonzeroload} that 
$a(t) \xrin k^* + \alpha/k_1$ for any $(a(0),b(0)) \in \R^2_{> 0}$. 
\end{proof}

\begin{theorem} \label{thm:flowmotif5exp_main}
Consider the mass action system \eqref{eq:motif5} with time-dependent inflow rates $g_a: \R_{\ge 0} \to \R_{\ge 0}$ and $g_b: \R_{\ge 0} \to \R_{\ge 0}$, such that 
\begin{equation} \label{eq:intcond5_main}
G(t) \coloneqq \int_{0}^{t} \left(g_a(s) + g_b(s)\right) ds  \xrin \infty. 
\end{equation}
We have that $a(t) \xrin k^*$ for any $(a(0),b(0)) \in \R^2_{> 0}$ if one of the following holds:
\been
\item $g_a$ is bounded, 
\item $g_a(t) \to \infty$ and 
$\ds \frac{\dot g_a(t)}{g_a(t)^{1/2} g(t)} \to 0$.
\enen
\end{theorem}
\begin{proof}
If $g_a$ is bounded, then $g_a(t)/(G(t))^2 \xrin 0$ and so the result follows immediately from Theorem \ref{thm:flowmotif5exp}. 
Suppose that $g_a(t) \xrin \infty$ and $\ds \frac{\dot g_a(t)}{g_a(t)^{1/2} g(t)} \to 0$. 
Since $g_a(t) \to \infty$, clearly $G(t) \to \infty$ and so by applying L'Hospital's rule, we have that 
\begin{align*}
 \lim_{t \to \infty} \frac{\sqrt{g_a(t)}}{G(t)}   = \lim_{t \to \infty} \frac{\frac{d}{dt} \sqrt{g_a(t)}}{\frac{d}{dt} G(t)} 
= \lim_{t \to \infty} \frac{\dot g_a(t)}{2g_a(t)^{1/2} g(t)} = 0, 
\end{align*}
and so $g_a(t)/(G(t))^2 \to 0$. Then from Theorem \ref{thm:flowmotif5exp}, $a(t) \xrin k^*$. 
\end{proof}

\begin{remark} \label{rem:flowmotif5exp_main}
When $g_a(t) \xrin \infty$, the limit of $\ds \frac{\dot g_a(t)}{g_a(t)^{1/2} g(t)}$ cannot be positive. Indeed suppose that  $\ds \frac{\dot g_a(t)}{g_a(t)^{1/2} g(t)} \xrin \alpha > 0$, then after some finite time $T_c > 0$, we have that for some positive constant $c > 0$, 
\[
\frac{d}{dt} g_a(t) > c \alpha g_a(t)^{1/4} g(t) \ge c \alpha g_a(t)^{5/4}, 
\]
and the differential equation $\dot g_a(t) = c \alpha g_a(t)^{5/4}$ has a finite time blow up. 
\end{remark}

\begin{remark}
A simpler result is obtained by replacing the term $\ds \frac{\dot g_a(t)}{g_a(t)^{1/2} g(t)}$ in Theorem \ref{thm:flowmotif5exp_main} and Remark \ref{rem:flowmotif5exp_main}  with the term  $\ds \frac{\dot g_a(t)}{g_a(t)^{3/2}}$. 
\end{remark}

\begin{example} \label{ex:tetrationgrowth}
Consider the following reaction network for some positive integer $m$ 
\begin{equation} \label{net:tetrationgrowth}
\begin{aligned} 
Z_1 &\xrightarrow{1} 2Z_1,\\
Z_{2} + Z_{1} &\xrightarrow{1} 2Z_{2} + Z_{1}, \\
Z_2 &\xrightarrow{1} Z_2 + A\\
A + 2B &\xrightarrow{k_1} 3B, \\
2B &\xrightarrow{k_2} A + B. 
\end{aligned}
\end{equation}
We assumed that the rate constants for some of the equations are $1$ for simplicity. 
The variable $z_1$ satisfies the differential equation $\dot z_1 = z_1$, which has an exponentially growing solution, $z_1 = e^t$, assuming initially $z_1(0) = 1$. 
Then $z_2$ satisfies the differential equation 
\[
\dot z_{2} = z_1 z_{2}, 
\]
which has the solution:
\[
z_{2}(t) =  e^{e^t},  
\]
assuming that $z_2(0) = e$. 
The ODE system satisfied by the species $A$ and $B$ is: 
\begin{equation}
\begin{aligned} 
\dot a &=  k_1b^2(k^*-a) + e^{e^t},\\
\dot b &= -k_1b^2(k^*-a). 
\end{aligned}
\end{equation}
Let $g_a(t) \coloneqq e^{e^t}$, then clearly $g_a(t) \to \infty$ and 
\[
\frac{\dot g_a(t)}{g_a(t)^{3/2}}  = \frac{e^t e^{e^t}}{e^{3e^t/2}} = e^{t-e^t/2} \xrin 0. 
\]
Then, by Theorem \ref{thm:flowmotif5exp_main}, $a \to k^*$. 
\end{example}

\subsection{Zero slope motif} \label{subsec:motif3}

We consider the motif ~
\begin{tikzpicture}[scale=0.5]
\draw [->, line width=1.25] ({-3+0.5},{0}) -- ({-3-0.5},{0});
\draw [->, line width=1.25] ({-3-0.5},{0}) -- ({-3+0.5},{0});
\draw [-, line width=1, dashed] ({-3-0.5},{0}) -- ({-3+0.5},{0});
\end{tikzpicture}
~
and the following open reaction network: 

\begin{equation}
\begin{aligned} 
B \xrightleftarrows{k_2}{k_1} A+B, \\
0    \xrightarrow{g_a(t)} A,   \quad &
0    \xrightarrow{g_{b}(t)} B. 
\end{aligned}
\end{equation}
whose mass action system is 
\begin{equation} \label{eq:motif3}
\begin{aligned} 
\dot a &=  k_1b(k^*-a) + g_a(t),\\
\dot b &=  g_b(t), 
\end{aligned}
\end{equation}
where $k^* = k_2/k_1$.

\begin{theorem}
Consider \eqref{eq:motif3} with $g_a, g_b:\R_{\ge 0} \to \R_{\ge 0}$ such that $\int_0^t b(s) ds \xrin \infty$. 
Then 
\[
a(t) \xrin k^* + \frac{1}{k_1} \lim_{t \to \infty} \frac{g_a(t)}{b(0) + \int_0^t g_b(s) ds}, 
\]
provided the limit exists as a finite value. 
\end{theorem}
\begin{proof}
The result follows immediately from applying Theorem \ref{thm:nonzeroload}. 
\end{proof}
\begin{remark}
It is worth contrasting the dynamic ACR properties of the negative slope motif with the other motifs. The negative slope motif has particularly striking robustness properties when inflows are added, which is not the case for the other motifs. 
This is further argument in favor of the attention that the negative slope motif has attracted, besides having a positive conservation law and being relatively simple in nature. 
\end{remark}

\subsection{Enzyme catalysis with bifunctional enzyme} \label{subsec:bifun}

Bifunctional enzymes are those which have two different (usually opposing) functions, for instance helping both to catalyze production as well as elimination of a substrate \cite{alon2019introduction,joshi2023bifunctional}. For instance, a bifunctional enzyme might facilitate phosphorylation and when in a different form, facilitate dephosphorylation. 
Bifunctional enzymes have been shown to play a critical role in generating static absolute concentration robustness \cite{hart2013utility,joshi2023bifunctional}. Here we show that an enzyme catalysis network with a bifunctional enzyme can produce dynamic absolute concentration robustness as well. 

We consider a mass action system with inflows that results from the reaction network depicted below. 
\begin{equation} \label{net:simp_futile}
\begin{aligned} 
0 \xrightarrow{g_x} X,~ 0 \xrightarrow{g_y} Y,~& 0 \xrightarrow{g_e} E,~ 0 \xrightarrow{g_c} C, \\
X+E  \stackrel[k_2]{k_1}{\rightleftarrows} &C \stackrel{k_3}{\to} Y+E, \\ 
Y+C & \xrightarrow{k_4} X+C. 
\end{aligned}
\end{equation}
The system of ODEs for the mass action system is
\begin{equation} \label{eq:simp_futile}
\begin{aligned} 
\dot x &= -k_1xe+k_2c+k_4cy +g_x, \\
\dot y &= k_3c-k_4cy +g_y, \\
\dot e &= -k_1xe+(k_2+k_3)c +g_e,  \\
\dot c &= k_1xe-(k_2+k_3)c + g_c. 
\end{aligned}
\end{equation}
We note that the variable $y$ is in power-engine-load form: 
\begin{align*}
\dot y 
= {\color{olive} k_4 c} {\color{magenta} \left( \frac{k_3}{k_4}- y\right)} + {\color{cyan} g_y}, 
\end{align*}
and therefore is a candidate for dynamic ACR. 
The key to proving dynamic ACR for different inflows or no inflows is to show that the ``power'' variable $c$ is eventually bounded away from $0$ for any initial value. 
We will do this here for a variety of cases with inflows. 
When there is no inflow in $Y$ ($g_y = 0$) and assuming that $Y$ is a dynamic ACR species, then it is easy to see that the dynamic ACR value must be $k^* \coloneqq k_3/k_4$, the same as its static ACR value. 
In the case of no inflows (in any species $g_x = g_y = g_e = g_c = 0$), one can further show that $y$ has dynamic ACR, that is for any positive initial value with $T_x \coloneqq x(0) + y(0) + c(0) \ge k^*$, we have that  $y \to k^*$. In other words, every initial value compatible with the hyperplane $\{y = k^*\}$ results in convergence to the hyperplane. 
We defer the proof of this result to future work because it requires additional ideas not directly related to the main theme of the paper. We focus here on certain interesting inflow cases.

In the theorem below, we show that $y$ can converge to the same value $k^*$ even in cases when some (or all) of the inflows are present. 

\begin{theorem}
Consider \eqref{eq:simp_futile} with $g_c > 0$ (while the other inflow rates $g_x, g_y, g_e$ are arbitrary nonnegative constants). Then $y \xrin k^* \coloneqq k_3/k_4$.  
\end{theorem}
\begin{proof}
Introduce the invertible change of variables $(x,y,e,c) \mapsto (y,c,T_x,T_e) \in \R^4_{\ge 0}$, where $T_x  \coloneqq  x+ y + c$ and $T_e  \coloneqq  e+c$. The new variables satisfy
\begin{equation} \label{eq:simp_futile_red}
\begin{aligned} 
\dot y &= k_4c\left(k^* - y \right) +g_y \\
\dot c &= k_1\left((T_x - y - c)(T_e - c)- k c \right) + g_c \\
\dot T_x &= g_x + g_y + g_c =: g_1 \\
\dot T_e &= g_e + g_c =: g_2
\end{aligned}
\end{equation}
where $k  \coloneqq  (k_2+k_3)/k_1$, $k^* \coloneqq k_3/k_4$, as well as $0 \le c(t) \le T_e(t)$ and $0 \le c(t) + y(t) \le T_x(t)$ for all $t \ge 0$. 
The hypothesis $g_c > 0$ implies that $g_1 >0$ and $g_2 > 0$. 
The evolution of $T_x$ and $T_e$ is uncoupled from that of $y$ and $c$, and is easily seen to have the solution $T_x(t) = T_x(0) + g_1 t$ and $T_e(t) = T_e(0) + g_2 t$. 
So it is natural to think of \eqref{eq:simp_futile_red} as a two-dimensional non-autonomous system in  the variables $(y,c)$. 

Since the first term of $\dot c$ is nonnegative, we have that $\dot c \ge -(k_2 + k_3) c + g_c$ which is positive for any $c < g_c/(k_2 + k_3)$ and so after some finite time, $c(t) > g_c/(k_2 + k_3)$. But this implies that $\dot y < 0$ for any sufficiently large $y$, and so $y$ is bounded above by some $\bar y > 0$ after some finite time. 

Now suppose that 
\[
c(t) < \frac12 \left( \min \left\{T_x(t)  - \bar y, T_e(t) \right\} - \frac{k}2 \right). 
\]
Then 
\begin{align*}
\dot c &= k_1\left( (T_x(t) - y - c)(T_e(t) - c)- k c\right) + g_c \\
& > k_1\left( (T_x(t) - \bar y - c)(T_e(t) - c)- k c\right) + g_c \\
&> k_1 \left( \left(c+ k/2 \right)^2 - k c \right) + g_c \\
&= k_1 \left(c^2+ k^2/4 \right) + g_c, 
\end{align*}
which implies that $c \to \infty$ since $T_x(t) - \bar y \to \infty$ and $T_e(t) \to \infty$. 

The result $y \to k^*$ then follows from applying Theorem \ref{thm:nonzeroload}. 
\end{proof}

We now present a case where specific inflows result in $y$ converging to a finite value, which is however not the same as $k^*$ (unless $g_y = 0$), the ACR value of $y$ when there are no inflows.

\begin{theorem}
Consider \eqref{eq:simp_futile} with $g_c = 0$, $g_e = 0$, $g_y \ge 0$ and $g_x + g_y > 0$. Then 
\[
y \xrin  k^* +  \frac{g_y}{k_4(e(0)+c(0))}.
\]  
In particular, $y \to k^*$ if $g_y = 0$.
\end{theorem}
\begin{proof}
Introduce the invertible change of variables $(x,y,e,c) \mapsto (y,c,T_x,T_e) \in \R^4_{\ge 0}$, where $T_x  \coloneqq  x+ y + c$ and $T_e  \coloneqq  e+c$. The new variables satisfy
\begin{equation} \label{eq:simp_futile_red2}
\begin{aligned} 
\dot y &= k_4c\left(k^* - y \right) +g_y \\
\dot c &= k_1\left((T_x(t) - y - c)(T_e - c)- k c \right)  \\
\dot T_x &= g_x + g_y =: g_1 > 0 \\
T_e(t) &= T_e(0) \coloneqq T_e
\end{aligned}
\end{equation}
where $k  \coloneqq  (k_2+k_3)/k_1$, $k^* \coloneqq k_3/k_4$, as well as $0 \le c(t) \le T_e$ and $0 \le c(t) + y(t) \le T_x(t)$ for all $t \ge 0$. 

Suppose first that $y \le k^*$ for all time. 
For arbitrary $\ve$, suppose that $c < T_e - \ve$. Then we have that 
\begin{align*}
\dot c &= k_1\left( (T_x(t) - y - c)(T_e - c)- k c\right)  \\
& > k_1\left( (T_x(t) - k^* - c)\ve - k c\right)  \\
& > k_1\left( (T_x(t) - k^* - T_e)\ve - k T_e \right). 
\end{align*}
Since $T_x(t) \uparrow \infty$, we have that the right hand side is arbitrarily large after some finite time, which implies that $c(t) \to T_e$. 

Now suppose that $y(t) > k^*$ for some $t$. 
It is clear that $\{y>k^*\}$ is invariant for any $g_y \ge 0$. 
We may therefore assume, without loss of generality, that $y(t) > k^*$ for all $t \ge 0$. 
For convenience we add $0 = 0 \cdot c =  (k_3 - k_4 k^*)c$ to $\dot x$ and rewrite the original system \eqref{eq:simp_futile} as
\begin{equation}
\begin{aligned} 
\dot x & = -k_1xe+(k_2+k_3)c-k_4c(k^*-y) +g_x,  \\
\dot y &= k_4c(k^* - y) +g_y, \\
\dot e &= -k_1xe+(k_2+k_3)c,  \\
\dot c &= k_1xe-(k_2+k_3)c. 
\end{aligned}
\end{equation}
Let $u(t) \coloneqq k_1xe-(k_2+k_3)c = k_1 (xe - kc), ~ v(t) \coloneqq -k_4c(k^*-y)$, 
and so the previous system in terms of $u$ and $v$ is:
\begin{align} 
\dot x  = -u(t)+v(t) +g_x, \quad
\dot y = -v(t) +g_y, \quad
\dot e = -u(t), \quad
\dot c = u(t). 
\end{align}
Note that $y(t) > k^*$ is equivalent to $v(t) > 0$. Furthermore, 
\begin{equation} \label{tempineq}
\begin{aligned}
\dot u(t) &= k_1(\dot x e + x \dot e - k \dot c) \\
&= k_1(-eu + ev + g_x e - xu - ku) = k_1( -u(e+x+k) + e(v+ g_x)) \\
&> - k_1(e+x+k) u
\end{aligned}
\end{equation}
This shows that either $u(t) > 0$ for all time after some finite time has passed, or $u(t) \uparrow 0$. 
We argue that in both cases $c$ must converge to a positive limit. 

First suppose that $u(t) \uparrow 0$. Then $\dot x \ge 0$ and $\dot e \ge 0$, and so $x(t) e(t)$ is nondecreasing for any $t \ge 0$. 
If $xe \uparrow \infty$, then $u(t) = k_1(xe-kc) \uparrow \infty$ since $c$ is bounded above. 
But this contradicts $u(t) = k_1(xe - kc) \uparrow 0$. 
So $xe$ is bounded above and therefore converges to a positive limit that is greater than $x(0) e(0)$. 
But then from the equation $\dot c = k_1xe-(k_2+k_3)c$, we have that $c(t)$ converges to a positive limit that is greater than $x(0) e(0)/k$. 

Now suppose that $u(t) > 0$ for all large enough $t$. 
Then $\dot c = u(t) > 0$ for large enough $t$, so that $c(t)$ is increasing. Since $c(t)$ is bounded above by $T_e$, $c(t)$ must converge to a positive limit. 

We conclude that $c$ converges to a positive limit in every case. 
From the equation $\dot y = k_4c(k^* - y) +g_y$, we conclude that $y$ converges to a positive limit. 
If $c$ converges to a limit less than $T_e$, then from $\dot c = k_1\left((T_x(t) - y - c)(T_e - c)- k c \right)$, $c$ increases without bound because $T_x$ increases without bound. This is a contradiction. Therefore $c \to T_e$. 
%
We may therefore conclude from Theorem \ref{thm:nonzeroload} that 
\[
y \xrin  k^* +  \frac{g_y}{k_4(e(0)+c(0))}.
\]  
\end{proof}

\section{Discussion}

The property of dynamic ACR depends on global dynamics and this makes it difficult to prove it even in simple systems. 
It is helpful to have either network conditions or ODE characterizations. 
Extending our previous work on network motifs which was only in two dimensions, we provide a power-engine-load form applicable in arbitrary dimensions. 
We give convergence results and apply these to certain network models in the chemostat setting, where inflows and sometimes outflows are present. 
Even when the overall trajectory does not converge, we show that the ACR variable converges to a definite finite quantity, the so-called ACR value. 
Moreover, in many cases this is the same value that the system would converge to when taken as a closed system, i.e. without any inflows and outflows. 

This work has implications for understanding the dynamics of systems which are far from equilibrium. 
We focus here on showing dynamic ACR in a few representative examples. In future work, we apply some of these results to analyze more complex and biochemically realistic systems.

\subsection*{Acknowledgments}
B.J. was supported by NSF grant DMS-2051498. 
G.C. was supported by NSF grant DMS-2051568 and by a Simons Foundation fellowship.
We thank the referee for careful reading and helpful comments. 
\bibliographystyle{unsrt}
\bibliography{acr}

\end{document}